%% file: ideal_eng_red.tex
\documentclass[11pt]{article}
\usepackage{latexsym,amsthm,amsfonts,epsfig,psfrag }
\begin{document}

\newtheorem{prop}{Proposition}
\newtheorem{lemma}{Lemma}
\newtheorem{theorem}{Theorem}
\renewcommand{\proofname}{Proof}

\def \H{{\mathbb H}}
\def \R{{\mathbb R}}
\def \E{{\mathbb E}}
\def \Z{{\mathbb Z}}

\begin{center}
{\LARGE  \bf
 On simple ideal hyperbolic Coxeter polytopes }

\vspace{12pt} 
{\large Anna~Felikson, Pavel~Tumarkin}
\end{center} 

\vspace{12pt} 


      
\section*{Introduction}

Let $\H^n$ be the $n$-dimensional hyperbolic space and 
let $P$ be a simple polytope in $\H^n$. 
$P$ is called a {\it Coxeter polytope } if 
all dihedral angles of $P$ are submultiples of $\pi$.

Hyperbolic Coxeter polytopes are not classified yet.
Examples of compact Coxeter hyperbolic polytopes are known up
to dimension  $n\le 8$ only, and examples of non-compact finite volume
Coxeter polytopes are known up to dimension
$n\le 19$~\cite{V5},~\cite{V4} and $n=21$~\cite{Borch}.
It is also known that hyperbolic spaces of high dimension contain
no finite volume Coxeter polytope.
The estimate for the highest possible dimension of a finite volume
Coxeter polytope is based on the following result of V.~V.~Nikulin.

Let $P$ be an $n$-dimensional simple polytope (where ``simple'' means that
any $k$-dimensional face of $P$ belongs to exactly $n-k$ facets),
and $\alpha_i$, $i=0,1,\dots,n-1$, be the number of its $i$-dimensional
faces ($i$-faces for short). For any face $f$ of $P$
denote by  $\alpha_i^f$ the number of its $i$-faces.
Denote by
$$
\alpha_k^{(i)}=\frac{1}{\alpha_k}\sum_{\mbox{\scriptsize dim} f=k} 
\alpha_i^f
$$
the average number of $i$-faces of a $k$-face of $P$.

\begin{prop}[Nikulin~\cite{Nikulin}]
For any simple convex compact polytope $P$ in $\R^n$ 
for any $i<k\le[n/2] $ the following estimate holds:
$$
\alpha_k^{(i)}<{{n-i}\choose{n-k}}\frac{{{[n/2]}\choose i}+{{[(n+1)/2]}\choose i}}{{{[n/2]}\choose k}+{{[(n+1)/2]}\choose k}}.
$$

\end{prop}

Using this estimate for $2$-faces ($i=0$ and $k=2$) and the fact that
any compact Coxeter polytope is simple, Vinberg~\cite{Vinb}  
proved that 
no compact Coxeter polytope exists in $\H^n$ for $n> 29$.  

In~\cite{Khov}, Khovanskij proved that Nikulin's estimate holds for edge-simple polytopes 
(a polytope is called {\it edge-simple} if any edge is the intersection of exactly $n-1$ facets).
This was used by Prokhorov~\cite{Prokh} when he proved that no Coxeter polytope of finite volume exists in $\H^n$
for $n\ge 996$.

A polytope $P$ is called {\it ideal } if all vertices of $P$ belong to the 
boundary of $\H^n$. 

In this paper, we study simple ideal hyperbolic Coxeter
polytopes.
The main result is the following theorem.

\begin{theorem}
No simple ideal Coxeter polytope exists in $\H^n$ when $n> 8$.

\end{theorem}
     
Section~\ref{prelim} contains basic definition and facts concerning
Coxeter diagrams of spherical, Euclidean, and hyperbolic Coxeter
polytopes. In Section~\ref{3-4}, we study the combinatorics of 
Coxeter diagrams of simple ideal hyperbolic Coxeter polytopes.
We show that if $n>5$ then such a polytope has no triangular
$2$-faces and a few quadrilateral $2$-faces.
As shown in section~\ref{res}, if $n>8$ this contradicts Niculin's
estimate. 

The paper was mainly worked out in
the Max-Planck Institute for Mathematics in Bonn.
The authors are grateful to the Institute for hospitality. 

\section{Coxeter diagrams}
\label{prelim}

It is convenient to describe Coxeter polytopes in terms of {\it Coxeter diagrams}.

A {\it Coxeter diagram} is one-dimensional simplicial complex
with weighted edges, where weights are either of the type
$\cos\frac{\pi}{m}$ for some integer $m\ge 3$ or positive real numbers
no less than one. We can suppress the weights but indicate the same information by
labeling the edges of a Coxeter diagram in the following way: if
the weight $w_{ij}$ equals $\cos\frac{\pi}{m}$, $v_i$ and $v_j$
are joined by an $(m-2)$-fold edge or a simple edge labeled by
$m$; if $w_{ij}=1$, $v_i$ and $v_j$ are joined by a bold edge; if
$w_{ij}>1$, $v_i$ and $v_j$ are joined by a dotted edge labeled
by its weight.

A {\it subdiagram} of a Coxeter diagram $\Sigma$ is a subcomplex with the same weights 
as in $\Sigma$. 

Let $\Sigma$ be a diagram with $d$ nodes $u_1$,...,$u_d$. Define a
symmetric $d\times d$ matrix  $G(\Sigma)$ in the following way: 
$g_{ii}=1$; if two nodes  $u_i$ and $u_j$ are joined by an edge with
weight $w_{ij}$ then $g_{ij}=-w_{ij}$; if two nodes  $u_i$ and
$u_j$ are not adjacent then $g_{ij}=0$. 

A {\it Coxeter diagram $\Sigma(P)$ of Coxeter polytope} $P$ is a Coxeter
diagram whose matrix $G(\Sigma)$ coincides with Gram matrix of $P$. In
other words, nodes of Coxeter diagram correspond to facets of $P$. Two nodes
are joined by either $(m-2)$-fold edge or  $m$-labeled edge if the 
corresponding dihedral angle equals $\frac{\pi}{m}$. If the
corresponding facets are parallel the nodes are joined by a bold edge,
and if they diverge then the nodes are joined by a dotted edge.


By the order of the diagram we mean the number of its nodes.
By signature and rank of diagram $\Sigma$ 
we mean the signature and the rank of the matrix
$G(\Sigma)$.

A Coxeter diagram $\Sigma$ is called {\it elliptic} if 
the matrix $G(\Sigma)$ is positive definite. 
A connected Coxeter diagram $\Sigma$ is called {\it parabolic} if 
the matrix $G(\Sigma)$ is degenerate, and any subdiagram of $\Sigma$ is
elliptic. Elliptic and connected parabolic diagrams are exactly Coxeter diagrams
of spherical and Euclidean Coxeter simplices respectively, they were
classified by Coxeter~\cite{Cox}. 
We represent the complete list of elliptic and connected parabolic diagrams in
Table~\ref{el-par}.

\input cox_e.txt  

A non-connected diagram is called {\it parabolic} 
if it is a disjoint union of 
connected parabolic diagrams. A diagram is  called {\it indefinite} 
if it contains 
at least one connected component that is neither elliptic nor parabolic.



Let $F$ be a $k$-dimensional face of $P$.
Since $P$ is simple, the face $F$ belongs to exactly $n-k$ facets
$f_1,\dots,f_{n-k}$. Denote by $v_1,\dots,v_{n-k}$ the corresponding 
nodes of $\Sigma(P)$. Let $\Sigma_F$ be a subdiagram of  $\Sigma(P)$ 
with nodes $v_1,\dots,v_{n-k}$.
We say that $\Sigma_F$ is the {\it diagram of the face} $F$. 
By {\it complete diagram of the face} $F$ we mean the minimal subdiagram
of $\Sigma(P)$ containing the diagrams of 
all vertices of $F$.

The following properties of $\Sigma(P)$ and $\Sigma_F$ are proved in~\cite{Vinb1}.

\begin{itemize}

\item[$\bullet$] [Cor. of Th. 2.1] the signature of $G(\Sigma(P))$
equals $(n,1)$;

\item[$\bullet$] [Cor. of Th. 3.1] if a $k$-face $F$ is not an ideal 
vertex of $P$ (i.e. $F$ is not 
a point at the boundary of $\H^n$), then   $\Sigma_F$ is an elliptic diagram of rank $n-k$;

\item[$\bullet$] [Cor. of Th. 3.2] if $F$ is an ideal vertex of $P$ 
then   $\Sigma_F$ is a parabolic 
diagram of rank $n-1$; if $F$ is  a simple ideal vertex of $P$ 
(i.e. $F$ belongs to exactly $n$ facets)
then $\Sigma_F$ is connected;

\item[$\bullet$] [Cor. of Th. 3.1 and Th. 3.2] any elliptic subdiagram
of $\Sigma(P)$ corresponds to a face of $P$; 
any parabolic subdiagram of $\Sigma(P)$ is a subdiagram of the diagram 
of exactly one ideal vertex of $P$.  

\end{itemize}

For a simple ideal Coxeter polytope $P\subset\H^n$ this implies that

\begin{itemize}

\item[{\sc (i)}\phantom {\sc i}] 
Any two non-intersecting indefinite subdiagrams of $\Sigma(P)$ 
are joined in $\Sigma(P)$.

\item[\,{\sc (ii)}\,] 
Any elliptic subdiagram of  $\Sigma(P)$ contains at most $n-1$ nodes.

\item[{\sc (iii)}] Any parabolic subdiagram of  $\Sigma(P)$ is connected and contains exactly $n$ nodes.

\end{itemize}


\begin{lemma}
\label{nok}
A Coxeter diagram of a simple ideal Coxeter polytope in $\H^n$, $n>3$, 
contains only 
simple edges, $2$-fold edges and dotted edges.

\end{lemma}

\begin{proof}
It follows from Table~\ref{el-par} that any connected parabolic diagram 
containing at least three nodes
contains neither bold edges nor edges of multiplicity $m>2$.
Thus, its enough to show that any non-dotted edge of $\Sigma(P)$
belongs to some connected parabolic subdiagram of order $n$.
Indeed, such an edge (denote it by $uv$) together with its ends
compose
a rank 2 elliptic subdiagram. Hence, it is a diagram of some
$(n-2)$-face $F$ of $P$. The face $F$ has at least one vertex,
and the diagram of this vertex is a connected parabolic subdiagram of 
$\Sigma(P)$ of order $n$ containing the diagram of $F$, i.e.
containing the edge $uv$.

\end{proof}

\subsubsection*{Notation}

Let $F$ be a $k$-face of $P$ and let $f_1,\dots,f_{n-k}$ be the facets 
of $P$ containing $F$.
Let $v_1,\dots,v_{n-k}$ be the corresponding nodes of $\Sigma(P)$.
As above, we denote by  $\Sigma_F$  the diagram of the face $F$,
i.e. the subdiagram of $\Sigma(P)$ spanned by the nodes
 $v_1,\dots,v_{n-k}$.

$\bullet$
We write $\Sigma_F=\;<\!\!v_1,\dots,v_{n-k}\!\!>$ and 
$\Sigma_F=\;<\!\!v_1,\Theta\!\!>$, 
where $\Theta=\;<\!\!v_2,\dots,v_{n-k}\!\!>$. 
We denote by $\Sigma\!\setminus\!\{v_1,...,v_m\}$ the 
subdiagram of $\Sigma$ spanned by all nodes of $\Sigma$ 
different from $v_1,...,v_m$.

$\bullet$
For elliptic and parabolic diagrams we use standard notation 
(see Table~\ref{el-par}).
For example, we write $\Sigma_F=\widetilde A_{n-1}$ if $F$ is an ideal
vertex of the type $\widetilde A_{n-1}$.

$\bullet$
Let  $v$ and $u$ be two nodes of $\Sigma(P)$.
We write 
\begin{itemize}

\item[]
$[v,u]=0$ if $u$ and $v$ are not joined in $\Sigma(P)$; 

\item[]
$[v,u]=1$ if $u$ and $v$ are joined by a simple edge; 

\item[]
$[v,u]=2$ if $u$ and $v$ are joined by a $2$-fold edge; 

\item[]
$[v,u]=\infty$ if $u$ and $v$ are joined by a dotted edge. 


\end{itemize}

\section{Absence of triangular 2-faces 
and estimate for quadrilateral 2-faces.}
\label{3-4}

Let $P$ be a simple ideal Coxeter polytope in $\H^n$ 
and let $V$ be a vertex of $P$. 
Since $P$ is simple, the vertex $V$ is contained in exactly $n$ edges 
$VV_i$, $i=1,\dots,n$.
Denote by $v_i$ the node of $\Sigma_V$ such that 
$\Sigma_{VV_i}=\Sigma_V\!\setminus\!\{v_i\}$.
Denote by $u_i$ the node of $\Sigma(P)$ such that $\Sigma_{V_i}=\;<\!\!u_i,\Sigma_{VV_i}\!\!>$.
Clearly, the diagram
$$\Sigma_i=<\!\! u_i, \Sigma_V\!\!> =<\!\! v_i, \Sigma_{V_i}\!\!> =
<\!\! u_i,v_i, \Sigma_{VV_i}\!\!>$$ is the complete diagram of the edge
$VV_i$.

Notice that $\Sigma_i$  contains exactly two parabolic subdiagrams
$\Sigma_V$ and $\Sigma_{V_i}$, therefore,
it is possible to find the nodes $v_i$ and $u_i$ in $\Sigma_i$ 
by formulae
$u_i=\Sigma_i\!\setminus\! \Sigma_V$, 
$v_i=\Sigma_i\!\setminus\! \Sigma_{V_i}$.
We say that a complete diagram of the edge  $\Sigma_i$ is 
{\it elementary} if there exists an automorphism of the diagram
$\Sigma_i$ interchanging the nodes
$u_i$ and $v_i$ and preserving the rest nodes. 
Otherwise we say that the complete diagram of the edge is 
{\it non-elementary}.

\begin{figure}[!h]
\begin{center}
\psfrag{a}{$v_i$}
\psfrag{b}{$u_i$}
\epsfig{file=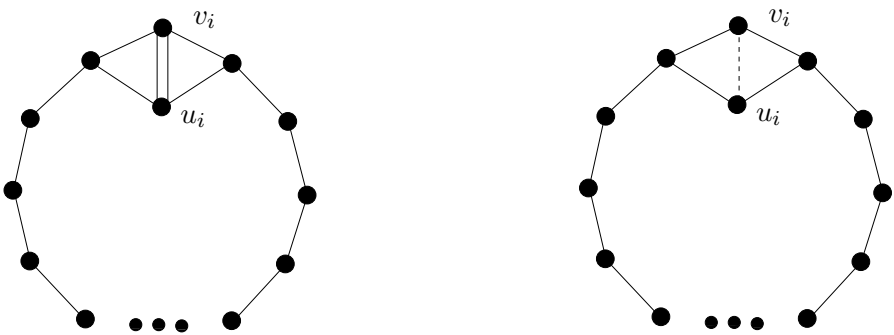,width=0.8\linewidth}
\end{center}
\caption{Two possibilities for the complete diagram
 $<\!\!v_i,u_i,\Sigma_{VV_i}\!\!>$ of the edge  $VV_i$, 
if $\Sigma_V=\widetilde A_{n-1}$, $n\ne 3, 8$,$9$. 
Both complete diagrams are elementary.
}
\label{A_n}
\end{figure}

For any connected parabolic diagram $\Sigma_V$ it is not difficult to 
describe all possible complete diagrams of edges containing $\Sigma_V$.
For example, suppose that $\Sigma_V=\widetilde A_{n-1}$, $n\ne 3,8,9$. 
Then $\Sigma_{VV_i}=\Sigma_V\!\setminus\! v_i=A_{n-1}$.
It is easy to see, that if $n\ne 3,8,9$ then $\widetilde A_{n-1}$ is 
the only parabolic diagram with $n$ nodes
containing a subdiagram $A_{n-1}$. Thus, $\Sigma_{V_i}=\widetilde
A_{n-1}$, and the diagram
$<\!\! v_i,u_i,\Sigma_{VV_i}\!\!>$ is elementary.  
Furthermore, $[v_i,u_i]\ne 0$ and  $[v_i,u_i]\ne 1$, otherwise
$<\!\!v_i,u_i,\Sigma_{VV_i}\!\!>$ does not satisfy condition {\sc (iii)}.
Hence, either $[v_i,u_i]=2$ or $[v_i,u_i]=\infty$ (Lemma~\ref{nok}), 
and the diagram $<\!\!v_i,u_i,\Sigma_{VV_i}\!\!>$ 
is one of two diagrams shown in Fig~\ref{A_n}.

Similarly, one can list all possible diagrams  
$<\!\!u_i,v_i,\Sigma_{VV_i}\!\!>$  for any other type of $\Sigma_V$ 
(recall that $\Sigma_V$ is one of the diagrams shown in the right
column of Table~\ref{el-par}).


\begin{lemma}
\label{non-elementary}
Let $<\!\!u_i,v_i,\Sigma_{VV_i}\!\!>$ be a non-elementary complete
diagram of the edge $VV_i$.
If $n>5$ and $\Sigma_{VV_i}$ is connected
then $<\!\!u_i,v_i,\Sigma_{VV_i}\!\!>$ is one of the diagrams
listed in Table~\ref{non-elem}.
 
\end{lemma}

\begin{table}[!h]
\begin{center}
\caption{Non-elementary diagram of the edge  
$<\!\!v_i,u_i,\Sigma_{VV_i}\!\!>$ such that $n>5$ and the diagram
$\Sigma_{VV_i}$ is connected. 
The waved edge connecting
$v_i$ and $u_i$ means that
$[v_i,u_i]\in \{ 0,1,2,\infty \}$
(for some of these values the conditions {\sc (i)}--{\sc (iii)} 
does not hold).}
\vspace{6pt}
\begin{tabular}{|c|c|c|}
\hline
&&\\
$\Sigma_{VV_i}=A_{n-1}$ & $\Sigma_{VV_i}=B_{n-1}$ &
$\Sigma_{VV_i}=D_{n-1}$
\\
&&\\
\hline
 
\begin{tabular}{c}
\epsfig{file=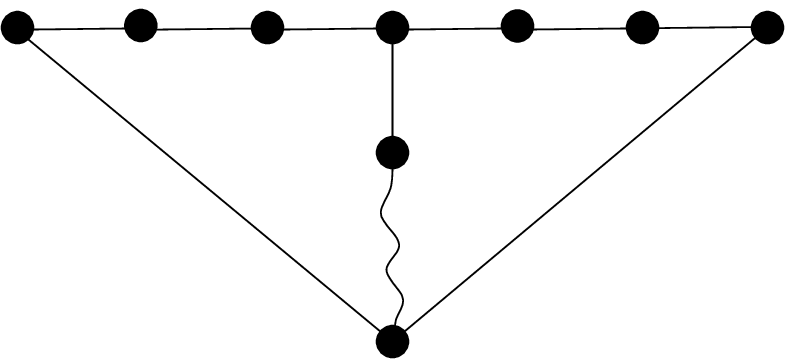,width=0.25\linewidth} \\
{\scriptsize 
$\widetilde A_7 \leftrightarrow \widetilde E_7$} \\ 
\\
\\
\epsfig{file=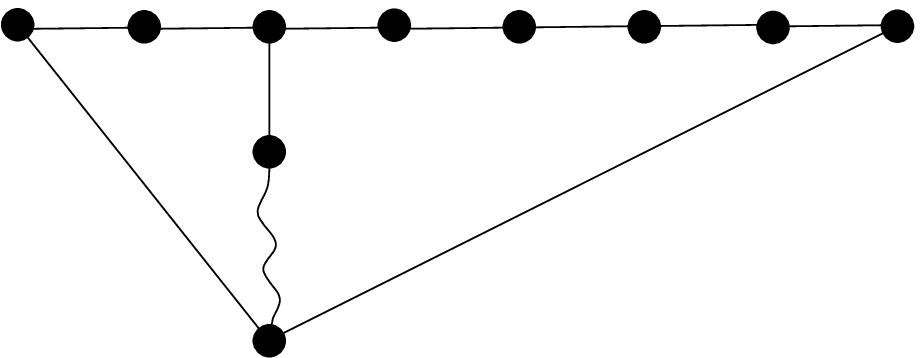,width=0.25\linewidth} \\
{\scriptsize 
$\widetilde A_8 \leftrightarrow \widetilde E_8$} \\ 
\\
\\
\epsfig{file=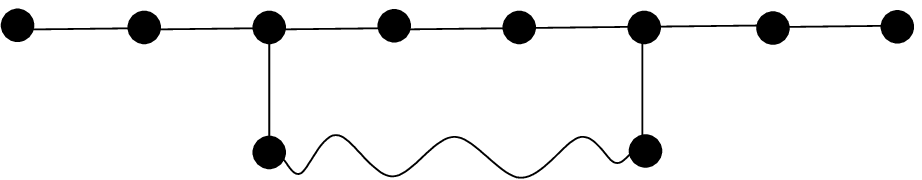,width=0.25\linewidth} \\
{\scriptsize 
$\widetilde E_8 \leftrightarrow \widetilde E_8$} \\

\end{tabular}
&
\begin{tabular}{c}
\epsfig{file=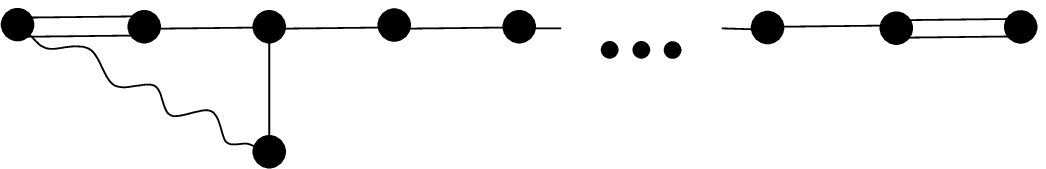,width=0.25\linewidth} \\
{\scriptsize 
$\widetilde C_{n-1} \leftrightarrow \widetilde B_{n-1}$} \\ 

\end{tabular}
&
 
\begin{tabular}{c}
\\
\epsfig{file=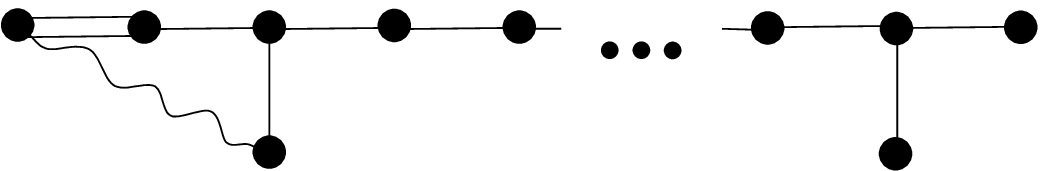,width=0.25\linewidth} \\
{\scriptsize 
$\widetilde D_{n-1} \leftrightarrow \widetilde B_{n-1}$} \\ 
\\
\epsfig{file=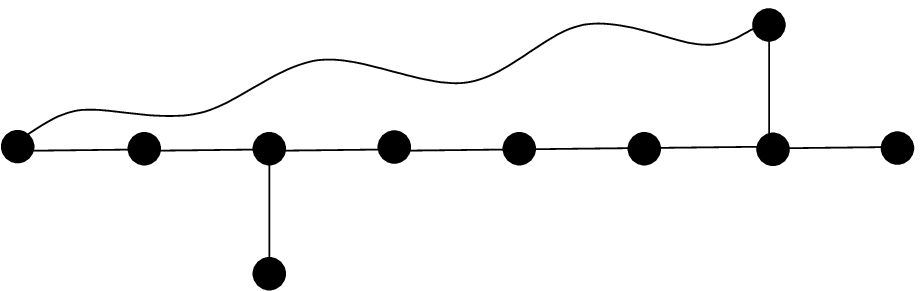,width=0.25\linewidth} \\
{\scriptsize 
$\widetilde D_8 \leftrightarrow \widetilde E_8$} \\ 
\\
\epsfig{file=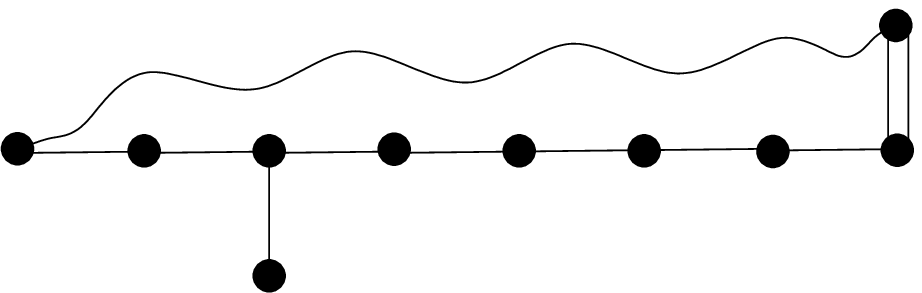,width=0.25\linewidth} \\
{\scriptsize 
$\widetilde B_8 \leftrightarrow \widetilde E_8$} \\ 
\\
\epsfig{file=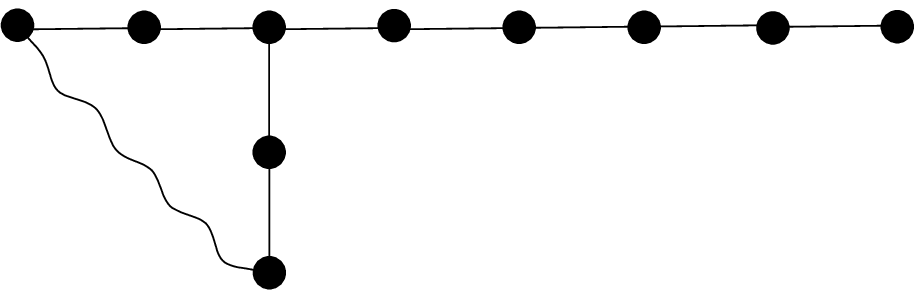,width=0.25\linewidth} \\
{\scriptsize 
$\widetilde E_8 \leftrightarrow \widetilde E_8$} \\ 

\end{tabular}
\\
&&\\
\hline

\end{tabular}
\label{non-elem}
\end{center}
\end{table}

\begin{proof}
There are two ways to obtain an edge  ${VV_i}$ with a non-elementary
complete diagram: either the diagrams of the vertices $V$ and $V_i$
are different or the diagrams are same but the nodes $u_i$ and $v_i$ 
are attached to the diagram $\Sigma_{VV_i}$ in different ways. 
Since $n>5$ and the diagram $\Sigma_{VV_i}$ is connected, 
the diagram $\Sigma_{VV_i}$ should be of one of the types
$A_{n-1}$, $B_{n-1}$, $D_{n-1}$ and $E_6, E_7, E_8$. 
Consider these cases.

{\bf 1.} Suppose that $\Sigma_{VV_i}=A_{n-1}$. Then the diagrams
$<\!\!u_i,\Sigma_{VV_i}\!\!>$ 
and $<\!\!v_i,\Sigma_{VV_i}\!\!>$ are parabolic diagrams of order $n$
containing a subdiagram of the type $A_{n-1}$.
Hence, each of these diagrams is of one of the types
$\widetilde A_{n-1}$ ($n\ge 6$),   $\widetilde E_7$ ($n=8$), 
and $\widetilde E_8$ ($n=9$).  
 
Furthermore, the diagram
$A_{n-1}$ extends to $\widetilde A_{n-1}$ in a unique way,
$A_{7}$ extends to $\widetilde E_7$ in a unique way, 
and $A_8$ extends to $\widetilde E_8$ in two different ways. 
Thus, if $n\ne 8,9$
the complete diagram of the edge $VV_i$ is always elementary.
If $n=8$ we obtain a unique non-elementary diagram
(where the multiplicity of the edge 
$u_iv_i$ may vary), denote this diagram by $\widetilde
A_7 \leftrightarrow \widetilde E_7$. If $n=9$ we obtain three
non-elementary diagrams 
(two diagrams of the type $\widetilde A_8 \leftrightarrow \widetilde
E_8$ and one of the type $\widetilde E_8 \leftrightarrow \widetilde E_8$), 
however, two diagrams of the type 
$\widetilde A_8 \leftrightarrow \widetilde E_8$ 
coincide  modulo the renumbering of the nodes.
So, in case  $\Sigma_{VV_i}=A_{n-1}$ we obtain three
non-elementary diagrams, see the left column of Table~\ref{non-elem}. 
 
{\bf 2.} Suppose that $\Sigma_{VV_i}=B_{n-1}$. 
Since $n>5$, each of the diagrams  
$<\!\!u_i,\Sigma_{VV_i}\!\!>$ and $<\!\!v_i,\Sigma_{VV_i}\!\!>$
is of the type $\widetilde B_{n-1}$ or  $\widetilde C_{n-1}$. 
The diagram $B_{n-1}$ may be extended to each of these diagram in a
unique way, and we obtain a unique non-elementary diagram
$\widetilde C_{n-1} \leftrightarrow \widetilde B_{n-1}$,
see the middle column of Table~\ref{non-elem}.

{\bf 3.} Suppose that $\Sigma_{VV_i}=D_{n-1}$. Then each of the diagrams  
$<\!\!u_i,\Sigma_{VV_i}\!\!>$ and $<\!\!v_i,\Sigma_{VV_i}\!\!>$
is of one of the types $\widetilde B_{n-1}$, 
$\widetilde D_{n-1}$ ($n\ge 6$), and
$\widetilde E_8$ ($n=9$). Since $n>5$, the diagram $D_{n-1}$
extends to each of the diagrams
$\widetilde B_{n-1}$ and $\widetilde D_{n-1}$ in a unique way.
The diagram $D_8$ extends to the diagram $\widetilde E_8$
in two different ways, and we obtain four non-elementary diagrams
$<\!\!u_i,v_i,\Sigma_{VV_i}\!\!>$ shown in the right column of
Table~\ref{non-elem}.

{\bf 4.} Suppose that $\Sigma_{VV_i}=E_6, E_7$ or $E_8$. 
Then each of the diagrams
$<\!\!u_i,\Sigma_{VV_i}\!\!>$ É $<\!\!v_i,\Sigma_{VV_i}\!\!>$
is of the types $\widetilde E_6$, $\widetilde E_7$ or $\widetilde E_8$
respectively. Since each of the diagrams $E_k$ ($k=6,7,8$)
extends to the diagram $\widetilde E_k$ in a unique  
(modulo the renumbering of the nodes) way,
the diagram $<\!\!u_i,v_i,\Sigma_{VV_i}\!\!>$ is elementary,
and the lemma is proved.

\end{proof}

The node $v$ of the diagram $\Sigma$ is called a {\it leaf} 
of the diagram  $\Sigma$,
if $v$ belongs to exactly one edge of $\Sigma$.

\begin{lemma}
\label{ab}
Let $<\!\!v_i,u_i,\Sigma_{VV_i}\!\!>$ be a complete diagram of the edge
$VV_i$.
If $n> 5$ then  $[v_i,u_i]\ne 0$ and $[v_i,u_i]\ne 1$.

\end{lemma}

\begin{proof}
Suppose that $[v_i,u_i]=0$ or $1$.

Assume that  the complete diagram $<\!\!v_i,u_i,\Sigma_{VV_i}\!\!>$ 
of the edge $VV_i$ is elementary, and consider two cases.

\begin{itemize}
\item[(a)] 
Suppose that $v_i$ is a leaf of  $<\!\!v_i,\Sigma_{VV_i}\!\!>$.
Denote by $a$ the node of $<\!\!v_i,\Sigma_{VV_i}\!\!>$
joined with $v_i$. Then the assumptions that the diagram 
 $<\!\!v_i,u_i,\Sigma_{VV_i}\!\!>$ is elementary and that
$[v_i,u_i]\ne\infty$ imply that
$<\!\!v_i,u_i,\Sigma_{VV_i}\!\setminus\! a\!\!>$ is an elliptic
subdiagram of order $n$, that contradicts condition
{\sc (ii)}.

\item[(b)] 
Suppose that $v_i$ is not a leaf of  $<\!\!v_i,\Sigma_{VV_i}\!\!>$.
Then there are at least two nodes  $a_1$ and  $a_2$ in 
$<\!\!v_i,\Sigma_{VV_i}\!\!>$ joined with $v_i$. 
Table~\ref{el-par} implies that one of the edges 
$a_1v_i$ and $a_2v_2$ is simple and another one is either simple or double.
Since the diagram  $<\!\!v_i,u_i,\Sigma_{VV_i}\!\!>$ is elementary and 
$[v_i,u_i]=0$ or $1$, we obtain that the diagram
$<\!\!v_i,u_i,a_1,a_2\!\!>$ contains a parabolic subdiagram of the type 
$\widetilde A_2$, $\widetilde C_2$ or $\widetilde A_3$,
which is impossible by condition {\sc (iii)}.

\end{itemize}

\noindent
Now, suppose that the diagram $<\!\!v_i,u_i,\Sigma_{VV_i}\!\!>$
is not elementary.
Suppose in addition that $\Sigma_{VV_i}$ is connected.
Then by Lemma~\ref{non-elementary}
the diagram $<\!\!u_i,v_i,\Sigma_{VV_i}\!\!>$ is one of the diagrams
listed in Table~\ref{non-elem}.
However, if $[v_i,u_i]=0$ or $1$, none of these diagrams satisfies  
conditions {\sc (ii)} and {\sc (iii)} simultaneously.

Therefore, the diagram $\Sigma_{VV_i}$ is not connected.
Let $\Sigma_1$ and $\Sigma_2$ be some connected components of 
$\Sigma_{VV_i}$
(it follows form Table~\ref{el-par} that  $\Sigma_{VV_i}$ contains at
most 3 connected components). Clearly, each of the nodes
$v_i$ and $u_i$ is joined with each connected component by exactly
one edge. Hence, the diagram $<\!\!\Sigma_1,\Sigma_2,v_i,u_i\!\!>$ 
contains a cycle $C$ including the nodes
$v_i$  and $u_i$. 

Suppose that the subdiagram
$<\!\!\Sigma_1,\Sigma_2\!\!>$ contains no double edges.
Then all edges of $<\!\!\Sigma_1,\Sigma_2,v_i,u_i\!\!>$ are simple,
and the cycle $C$ is a parabolic diagram of the type
$\widetilde A_k$ containing the nodes $u_i$ and $v_i$.
If $k<n-1$ this is impossible by condition  {\sc (iii)}, and the case
$k=n-1$ contradicts the assumption that  
$<\!\!v_i,u_i,\Sigma_{VV_i}\!\!>$ is the complete diagram of
the edge $VV_i$.

Therefore, at least one of the diagrams $\Sigma_1$ and $\Sigma_2$
contains a double edge which is included in the cycle $C$, 
i.e. either  $\Sigma_1$ or $\Sigma_2$ is a diagram
 $B_k$ for some $2\le k <n-1$.
We assume that $\Sigma_1=B_k$ and denote by
$t_1$ and $t_2$ the ends of the double edge in such a way that
$t_1$ is a leaf of $\Sigma_1$. Since the edge $t_1t_2$ belongs to the 
cycle $C$, one of the nodes
$u_i$ and $v_i$ (say, $u_i$) is joined with
$t_1$, and another one ($v_i$) is not. If $k>2$ then the nodes
 $t_1$ and $t_2$ are not  leaves of the parabolic diagram
$<\!\!u_i,\Sigma_{VV_i}\!\!>$, and hence, 
$<\!\!u_i,\Sigma_{VV_i}\!\!>=\widetilde F_4$ (see Table~\ref{el-par}),
that contradicts the assumption that $n>5$.  

Thus, $\Sigma_1=B_2=t_1t_2$, and $u_i$ is joined with $t_1$,
while  $v_i$ is joined with $t_2$. 
It follows from the classification of parabolic diagrams that the edges
$u_it_1$ and $v_it_2$ are simple. Consider two cases:
$[u_i,v_i]=1$ or  $[u_i,v_i]=0$.

\begin{itemize}
\item
If $[u_i,v_i]=1$ then the diagram 
$<\!\! u_i,v_i, \Sigma_2\!\!>$ also contains a cycle,
and by the same reasoning as above we obtain $\Sigma_2=B_2$.
Since the diagram $\widetilde C_{n-1}$ is the only connected
parabolic linear diagram of order $n$ containing a subdiagram
of the type $B_2+B_2$,  
the diagram $\Sigma_{VV_i}$ contains no other connected
components besides $\Sigma_1=B_2$ and  $\Sigma_2=B_2$. 
Therefore, $n-1=4$, which contradicts the assumption that $n>5$.

\item
If $[u_i,v_i]=0$ then $<\!\! u_i,v_i,\Sigma_1\!\!>=F_4$.
Suppose that $u_i$ and $v_i$ are joined with one and the same node $x$ 
of $\Sigma_2$. Then $<\!\! u_i,v_i,\Sigma_1,\Sigma_2\!\setminus\! x\!\!>$ 
is an elliptic diagram of  order $n$, which is impossible by condition  
 {\sc (ii)}.
Hence, $u_i$ and $v_i$ are joined with distinct nodes $x_1$ and $x_2$
of $\Sigma_2$. If $u_i$ is joined with $x_1$ by a simple edge, then
 $<\!\! x_1,u_i,t_1,t_2,v_i\!\!>=\widetilde F_4$ which contradicts 
either condition {\sc (iii)} or the assumption $n>5$. If $u_i$
is joined with $x_1$ by a double edge then $<\!\!
x_1,u_i,t_1,t_2\!\!>=\widetilde C_3$, which contradicts 
either condition {\sc (iii)} or the assumption $n>5$ again.

\end{itemize}

\end{proof}

\begin{lemma}
\label{triangles}
Let $P$ be a simple ideal Coxeter polytope in $\H^n$, $n> 5$.
Then $P$ has no triangular $2$-faces. 

\end{lemma}

\begin{proof}
Suppose that $UVW$ is a triangular $2$-face of $P$.
Since $P$ is simple, the triangle $UVW$ is contained in exactly 
$n-2$ facets.
There exists a unique facet containing the edge $VW$ and not
containing the triangle $UVW$. Denote this facet by $\bar u$.
Similarly, determine facets  $\bar v$ and $\bar w$
as facets containing the edges  $UW$ and $UV$ and not containing $UVW$.
Denote by $u$, $v$ and $w$ the nodes  of $\Sigma(P)$ corresponding 
to  $\bar u, \bar v$ and $\bar w$ respectively.
Notice that the diagram of $U$  coincides with
$<\!\!v,w,\Sigma_{UVW}\!\!>$. Similarly,
$\Sigma_V=\;<\!\!u,w,\Sigma_{UVW}\!\!>$ and
$\Sigma_W=\;<\!\!u,v,\Sigma_{UVW}\!\!>$ (see Figure~\ref{notation}a).
In particular, {\sc (iii)} implies that all these diagrams are
connected and parabolic. We also obtain that the diagram
$\Sigma=<\!\!u,v,w,\Sigma_{UVW}\!\!>$ is the complete diagram of the 
triangular $2$-face $UVW$ as well as the complete diagram of each of the
edges $VW$, $UW$ and $UV$.

\begin{figure}[!h]
\begin{center}
\psfrag{a}{$\bar u$}
\psfrag{b}{$\bar v$}
\psfrag{u1}{$\bar u_i$}
\psfrag{v1}{$\bar v_i$}
\psfrag{u2}{$\bar u_j$}
\psfrag{v2}{$\bar v_j$}
\psfrag{V}{$V$}
\psfrag{V1}{$V_j$}
\psfrag{V2}{$V_i$}
\psfrag{V12}{$V_{ij}$}
\psfrag{(a)}{a)}
\psfrag{(b)}{b)}
\psfrag{c}{$\bar w$}
\psfrag{q}{$Q$}
\psfrag{t}{$$}
\psfrag{A}{$U$}
\psfrag{A'}{$A'$}
\psfrag{C}{$W$}
\psfrag{B'}{$B'$}
\psfrag{B}{$V$}
\epsfig{file=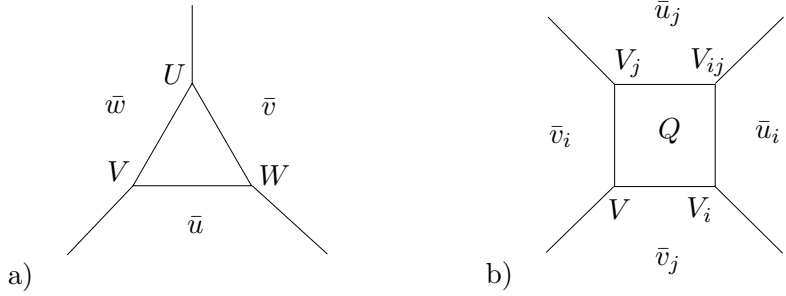,width=0.7\linewidth}
\end{center}
\caption{Notation  for a triangle (a) and for a quadrilateral (b). }
\label{notation}
\end{figure}

Consider the edge of $\Sigma$ joining $u$ and $v$.
By Lemma~\ref{ab}, either $[u,v]=2$ or $[u,v]=\infty$.

Suppose that $[u,v]=\infty$, then 
$\Sigma_W=\;<\!\!u,v,\Sigma_{UVW}\!\!>$ contains a dotted edge 
in contradiction to the assumption that $\Sigma_W$ is parabolic. 
Thus,  $[u,v]\ne \infty$, i.e. $[u,v]= 2$.
Similarly, $[v,w]= 2 $ and  $[u,w]= 2$.
Furthermore, since
$\Sigma_W$ is a parabolic diagram of order $n>5$, one of the nodes $u$ and
$v$ of the double edge $uv$ is a leaf. Assume that
$u$ is a leaf of $\Sigma_W$, i.e.  $u$ is not joined with
$\Sigma_{UVW}$. Then, evidently,  $v$ is joined with
$\Sigma_{UVW}$. Similarly, from the diagram $\Sigma_U=<\!\!
v,w,\Sigma_{UVW}\!\!>$ we obtain that $w$ is not joined with
$\Sigma_{UVW}$. Hence, the diagram
$<\!\!u,w,\Sigma_{UVW}\!\!>=\Sigma_{V}$ is not connected in
contradiction to condition {\sc (iii)}. 

\end{proof}

\begin{figure}[!h]
\begin{center}
\epsfig{file=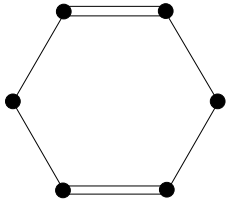,width=0.2\linewidth}
\end{center}
\caption{This diagram determines a 5-dimensional ideal Coxeter simplex.}
\label{5-simplex}
\end{figure}

Notice,  that an ideal Coxeter polytope in $\H^5$ may have a 
triangular $2$-face.
For example, the Coxeter diagram shown in Figure~\ref{5-simplex} determines 
a 5-dimensional ideal Coxeter simplex. All $2$-faces of any simplex 
are triangles.

\begin{lemma}
\label{quadrilaterals}
Let $V$ be a vertex of simple ideal Coxeter polytope $P$ in $\H^n$, $n>9$.
Then $V$ belongs to at most $n+3$ quadrilateral $2$-faces.

\end{lemma}

\begin{proof}
Let $Q$ be a quadrilateral $2$-face with vertices $V, V_i, V_j$ and $V_{ij}$. 
Then $Q$ belongs to $n-2$ facets, each edge of  $Q$ belongs to 
$n-1$ facets and each vertex belongs to $n$ facets.
Denote by $\bar v_i,\bar u_i,\bar v_j$ and $\bar u_j$ the facets 
not containing $Q$ and containing 
the edges $VV_j, V_iV_{ij},VV_i$ and $V_jV_{ij}$ respectively 
(see Figure~\ref{notation}b).
 Denote by $v_i,u_i,v_j$ and $u_j$ the nodes of $\Sigma(P)$ corresponding to the facets 
$\bar v_i,\bar u_i,\bar v_j$ and $\bar u_j$  respectively.

Then
$\Sigma_V=\;<\!\!v_i,v_j,\Sigma_Q\!\!>,\
\Sigma_{V_i}=\;<\!\!v_j,u_i,\Sigma_Q\!\!>,\
\Sigma_{V_j}=\;<\!\!v_i,u_j,\Sigma_Q\!\!>$, 
and $\Sigma_{V_{ij}}=\;<\!\!u_i,u_j,\Sigma_Q\!\!>$. 
Thus,
$$ \Sigma=\;<\!\!v_i,u_i,v_j,u_j,\Sigma_Q\!\!>$$ 
is the complete diagram of the face $Q$
(see Fig.~\ref{quad} for the example of a complete diagram
of  a quadrilateral).

\begin{figure}[!h]
\begin{center}
\psfrag{a}{$v_i$}
\psfrag{b}{$v_j$}
\psfrag{a'}{$u_i$}
\psfrag{b'}{$u_j$}
\psfrag{A}{$\Sigma_V$}
\psfrag{A'}{$<\!\!v_i,\Sigma_{V_i}\!\!>$}
\psfrag{C}{$<\!\!v_i,v_j,\Sigma_{V_{ij}}\!\!>$}
\psfrag{B'}{$<\!\!v_j,\Sigma_{V_j}\!\!>$}
\epsfig{file=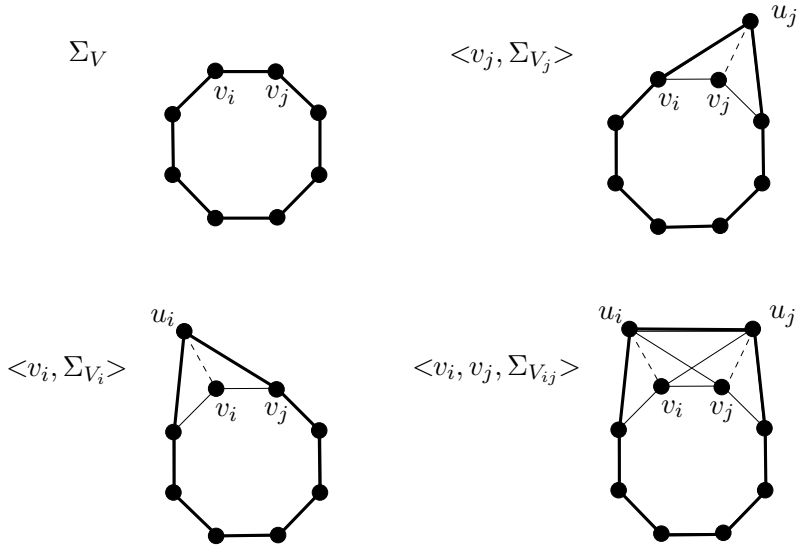,width=0.7\linewidth}
\end{center}
\caption{Example of a quadrilateral  $VV_iV_{ij}V_j$.}
\label{quad}
\end{figure}

Suppose that $\Sigma_V=\widetilde A_{n-1}$.
Since $n>8$, the diagram of each of the vertices 
$V_i,V_j,V_{ij}$ is of the type $\widetilde A_{n-1}$.
Consider the diagram $<\!\!v_i,u_i,\Sigma_{VV_i}\!\!>$, i.e. the
complete diagram of the edge  $VV_i$.
By Lemma~\ref{ab}, $[v_i,u_i]=\infty$ or $2$
(compare with Fig.~\ref{A_n}). Similarly, from the complete diagram of
the edge $VV_j$, we obtain that $[v_j,u_j]=\infty$ or $2$.

Suppose that the nodes $v_i$ and $v_j$ are not joined in $\Sigma_V$.
Then the diagram $\Sigma_Q=\Sigma_{V}\!\setminus\!\{v_i,v_j\}$ 
is not connected.
On the other hand,
$\Sigma_{V_{ij}}\!\setminus\!\{u_i,u_j\}=\Sigma_{Q}$, and we obtain
$[u_i,u_j]=0$ (see Fig.~\ref{aaa}). 
Since $n>5$, at least one of the connected components $\Sigma_1$ and
$\Sigma_2$ of the diagram $\Sigma_Q$ contains at least three nodes,
and we may assume that  $\Sigma_1=A_{k}$, $k\ge 3$. Denote by $w_i$
and $w_j$ the leaves of the diagram $\Sigma_1$ joining $\Sigma_V$ with 
$v_i$ and $v_j$ respectively. Then $\Sigma(P)$ contains two unjoined 
indefinite subdiagrams $v_iu_iw_i$ and $v_ju_jw_j$,
which is impossible by condition {\sc (i)}.
Therefore, the nodes $v_i$ and $v_j$ are joined in the diagram $\Sigma_V$,
that implies that each quadrilateral face containing the vertex
$V$ corresponds to a pair of  neighboring nodes in $\Sigma_V$.
Hence, $V$ belongs to at most $n$ quadrilateral $2$-faces.

\begin{figure}[!h]
\begin{center}
\psfrag{a}{$v_i$}
\psfrag{a'}{$u_i$}
\psfrag{a''}{$w_i$}
\psfrag{b}{$v_j$}
\psfrag{b'}{$u_j$}
\psfrag{b''}{$w_j$}
\psfrag{k}{$k_i$}
\psfrag{l}{$k_j$}
\epsfig{file=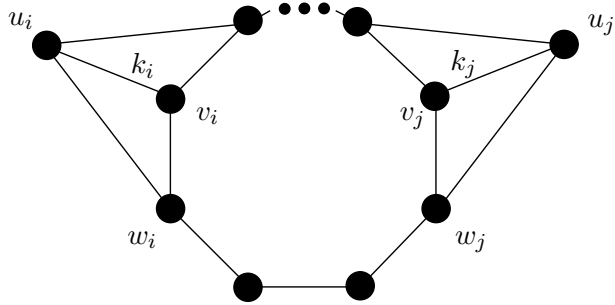,width=0.55\linewidth}
\caption{$v_iu_iw_i$ and $v_ju_jw_j$ are unjoined indefinite subdiagrams. In this diagram $k_i, k_j=2$ or $\infty$.} 
\label{aaa}
\end{center}
\end{figure}

From now on we assume that $\Sigma_V\ne \widetilde A_{n-1}$. Since $n>9$, 
$\Sigma_V=\widetilde B_{n-1}, \widetilde C_{n-1}$  or $\widetilde D_{n-1}.$
Define a {\it distance} $\rho(u,w)$ between two nodes $u$ and $w$ 
of connected 
graph as the number of edges in the shortest path connecting $u$ and $w$.

Let $x$ be a leaf of $\Sigma_V$. Denote by $\Sigma^{(5)}_{V}(x)$ 
a connected subdiagram of $\Sigma_V$ spanned by five nodes 
closest to the leaf $x$ in $\Sigma_V$ 
(i.e., if $v_k\in \Sigma^{(5)}_{V}(x)$ and  $v_l\notin \Sigma^{(5)}_{V}(x)$ 
then $\rho(x,v_k)\le \rho(x,v_l)$). Notice that if 
$\Sigma_V= \widetilde B_{n-1},\widetilde C_{n-1}$  or 
$\widetilde D_{n-1}$ when 
$n\ge 9$, then diagram $\Sigma^{(5)}_{V}(x)$ is well-defined for any
leaf $x$ of $\Sigma_V$.

Denote by $L(\Sigma_V)$ the set of leaves of $\Sigma_V$. 
Define 
$$\Sigma^{(5)}_{V}\stackrel{\mbox{\scriptsize def}}{=\!=}
\bigcup \limits_{x \in L(\Sigma_V)}\Sigma^{(5)}_{V}(x)$$ 
(see Fig.~\ref{end} for the example).
It is easy to see that  if $n> 10$ then $\Sigma^{(5)}_{V}$ consists of 
two connected components. 
If $n=10$, $\Sigma^{(5)}_{V}$ is connected. However, 
it contains two leaves $x$ and $y$ such that 
$\Sigma_V=\Sigma^{(5)}_{V}(x)\cup \Sigma^{(5)}_{V}(y)$,
and $\Sigma^{(5)}_{V}(x)\cap \Sigma^{(5)}_{V}(y)=\emptyset$.
In this case we say that the diagrams
$\Sigma^{(5)}_{V}(x)$ and $\Sigma^{(5)}_{V}(y)$ are ``components'',
and use this notion in case  $n=10$ instead of the connected
components in general case $n>10$.

\begin{figure}[!h]
\begin{center}
\epsfig{file=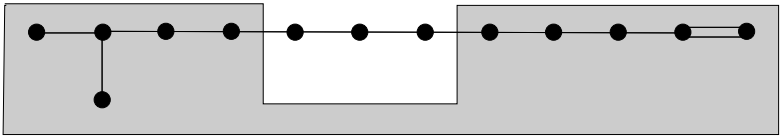,width=0.7\linewidth}
\end{center}
\caption{Subdiagram  $\Sigma^{(5)}_{V}$ for $\Sigma_V=\widetilde B_{12}$.}
\label{end}
\end{figure}

Suppose that $v_i$ and $v_j$ do not belong to the same connected 
component of  $\Sigma^{(5)}_{V}$ (respectively, to a ``component'' for  
$n=10$).
Suppose that  $v_i$ is not joined with $v_j$.
A direct check of 
the conditions {\sc (i)}--{\sc (iii)} for each of the possible diagrams
shows that if
$\Sigma_Q$ is a diagram of a quadrilateral $2$-face
then the corresponding connected component (or the ``component'')
of the diagram $\Sigma^{(5)}_{V}$ coincides (modulo interchanging
of $v_i$ and $v_j$) with one of the following diagrams:

\begin{center}
\psfrag{a}{$v_i$}
\psfrag{b}{$v_j$}
\epsfig{file=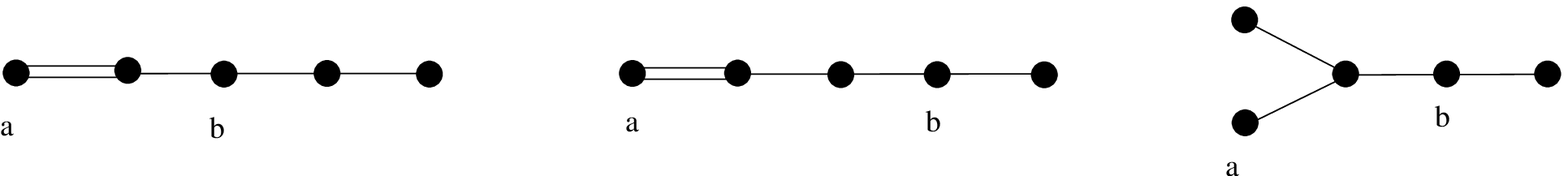,width=0.9\linewidth}
\end{center}

Therefore, a quadrilateral $2$-face $Q$ containing $V$
is of one of the following types:
either $Q$ corresponds to a pair of joined nodes in $\Sigma_V$ 
(there are $n-1$ of such pairs)
or to one of two pairs described above for each of the connected
components (or ``components'' for  $n=10$) 
of the diagram $\Sigma^{(5)}_{V}$. 
Hence, $V$ belongs to at most $2+2+(n-1)=n+3$ quadrilaterals.

\end{proof}

\begin{lemma}
\label{quad9}
Let $V$ be a vertex of a simple ideal Coxeter polytope $P$ in $\H^9$.
Then $V$ belongs to at most $15$ quadrilateral $2$-faces.

\end{lemma}

\begin{proof}

Let $v_1,\dots,v_9$ be the nodes of $\Sigma_V$.
While proving Lemma~\ref{quadrilaterals} we estimated the number
of pairs $(v_i,v_j)$ such that the $2$-face corresponding to the
diagram $\Sigma_{V}\!\setminus\! \{v_i,v_j\}$ {\it may be} 
quadrilateral. In other words, were looking for the pairs 
$(v_i,v_j)$, such that the diagram 
$\Sigma_V$ can be accompanied by some additional nodes $x$ and $y$ 
subject to the following two properties:
1) the diagram $<\!\!x,y,\Sigma_{V}\!\!>$ 
satisfies to conditions {\sc (i)---(iii)}, 2) the diagrams
$<\!\!x,\Sigma_{V}\!\setminus\! v_i\!\!>$, $<\!\!y,\Sigma_{V}\!\setminus\!
v_j\!\!>$ and $<\!\!x,y,\Sigma_{V}\!\setminus\! \{v_i,v_j\}\!\!>$
are connected and parabolic.
In this case 
$<\!\!x,y,\Sigma_{V}\!\!>$  is a complete diagram of a quadrilateral
$2$-face with diagram
$\Sigma_{V}\!\setminus\! \{v_i,v_j\}$.

Clearly, this estimate of the number of quadrilaterals is rough.
In particular, for  $n=9$ the result of this estimate worse
than one claimed in the lemma. 
To prove the lemma we use another method leading to the better
estimate, but using much more computations.
We proceed by the following algorithm:
 
\begin{itemize}

\item[Step 1.]
We consider the cases $\Sigma_V=\widetilde A_8,\widetilde B_8,
\widetilde C_8,\widetilde D_8$ and $\widetilde E_8$ separately.

\item[Step 2.]
We want to list all possibilities for 
complete diagrams of edges incident to $V$. 
To do this, for each node
$v_i$ of $\Sigma_V$ ($i=1,\dots,9$) 
we list all possible diagrams  $<\!\!u_i^k,\Sigma_{V}\!\!>$,
$k=1,\dots,k_i$ satisfying conditions {\sc (i)---(iii)} and such that
$<\!\!u_i^k,\Sigma_{V}\!\setminus\! v_i\!\!>$ is a parabolic diagram.
Here $k_i$ stays for the number of different complete diagrams of
edges found for each of the nodes $v_i$ of $\Sigma_V$. 
A straightforward check shows that
 $1\le k_i\le 8$ for different nodes of the diagrams
$\widetilde A_8,\widetilde B_8,\widetilde C_8,\widetilde D_8$ and
$\widetilde E_8$ (for instance, $k_i=8$ for one of the nodes of $\widetilde
E_8$). 

\item[Step 3.]
For each pair $(<\!\!u_i^k,\Sigma_{V}\!\!>$,
$<\!\!u_j^l,\Sigma_{V}\!\!>)$ of complete diagrams of edges
($i\ne j$) we check if it is possible to assign a weight to the edge
$u_i^k u_j^l$ in order to turn the diagram 
$<\!\!u_i^k,u_j^l,\Sigma_{V}\!\!>$ into a complete diagram of a 
quadrilateral.  In particular, this implies that the diagram
$<\!\!u_i^k,u_j^l,\Sigma_{V}\!\setminus\!\{v_i,v_j\}\!\!>$ is
parabolic, and hence, the nodes $u_i^k$ and $u_j^l$ are either 
unjoined or joined by a simple or double edge.
Each of the pairs  $(<\!\!u_i^k,\Sigma_{V}\!\!>$,
$<\!\!u_j^l,\Sigma_{V}\!\!>)$ obtained we call a {\it good} pair of
complete diagrams of edges.

\item[Step 4.]
For each of the nodes $v_i$, $1\le i\le 9$, we choose a complete
 diagram $<\!\!u_i^{r_i},\Sigma_{V}\!\!>$ of an edge, $1\le r_i\le k_i$. 
Then compute the total number of the good pairs 
$(<\!\!u_i^{r_i},\Sigma_{V}\!\!>,<\!\!u_j^{r_j},\Sigma_{V}\!\!>)$ 
 (where $i\ne j$). At this step we should check rather huge number of
cases (more than 15000 in case $\Sigma_{V}=\widetilde E_8$), 
therefore, this was done by a computer program.

The number obtained in this step we denote by $M(r_1,\dots,r_9)$.  
Denote by
$M(\Sigma_V)$ the maximal value of $M(r_1,..,r_9)$ on the 9-tuples
$(r_1,..,r_9)$, where $1\le r_i\le k_i$. Clearly,
the number of quadrilateral $2$-faces containing the vertex $V$
is bounded by $M(\Sigma_V)$. Notice also, that the estimate is still
rough (for example, we do not check if the conditions  
{\sc (i)---(iii)} are satisfied by subdiagrams containing more than
$n+2$ nodes). The computation shows that
\begin{center}
$M(\widetilde A_8)=15$,\\
$M(\widetilde B_8)=14$,\\
$M(\widetilde C_8)=12$,\\
$M(\widetilde D_8)=15$,\\
$M(\widetilde E_8)=14$.\\
\end{center}

\end{itemize}

Thus, for any type of $\Sigma_V$ we obtain that
$V$ belongs to at most $15$ quadrilateral $2$-facets.

\end{proof}

\section{Absence of simple ideal Coxeter polytopes \\ in large dimensions. }
\label{res}

Recall that $\alpha_i$ denotes the number of $i$-faces of a polytope
$P$ and $\alpha_k^{(i)}$ denotes the average number of $i$-faces of  
$k$-face of $P$.


\begin{lemma}
\label{l}
Let $P$ be an $n$-dimensional simple polytope and let $l$ be the
number of vertices of $P$. Then 
$$
\frac{l}{\alpha_2}=\frac{2}{n(n-1)}\alpha_2^{(1)}.\eqno (1)
$$

\end{lemma}

\begin{proof}
Denote by $m_i$ the number of $i$-angular $2$-faces of $P$.
Let us compute the total number $N$ of vertices of $2$-faces.
Clearly, $N=\sum \limits_{i\ge 3} i\cdot m_i $.
On the other hand,
each pair of edges incident to one vertex of  simple polytope determines a 
$2$-face of the polytope. Thus, $N=l\frac{n(n-1)}{2}$, and we obtain
the following equality
$$l\frac{n(n-1)}{2}=\sum \limits_{i\ge 3}{ i\cdot m_i }.\eqno (2)$$
By definition, 
$$
 \alpha_2^{(1)}=\frac{\sum \limits_{i\ge 3} i\cdot m_i }{\alpha_2}.\eqno (3)
$$
Combining (2) and (3), we obtain 
$$ 
\frac{l}{\alpha_2}=\frac{2}{n(n-1)}\frac{\sum \limits_{i\ge 3} i\cdot m_i}{\alpha_2}= 
\frac{2}{n(n-1)}\alpha_2^{(1)}.
$$
\end{proof}

\noindent
{\bf Proof of the theorem.\ }
We use the notation from Lemma~\ref{l}. Recall, that 
$\alpha_2=\sum \limits_{i\ge 3} m_i$. 
By Lemma~\ref{triangles}, $m_3=0$. Using (3), we obtain
$$
\alpha_2^{(1)}\ge \frac{1}{\alpha_2}(4m_4+5\sum\limits_{i\ge 5} m_i)=
 \frac{1}{\alpha_2}(5\sum\limits_{i\ge 4}m_i - m_4)=
5-\frac{m_4}{\alpha_2}.
\eqno(4)
$$

Consider Nikulin's estimate for $\alpha_2^{(1)}$:
$$
%
\alpha_2^{(1)}<{{n-1}\choose {n-2}}\frac{{{[n/2]}\choose 1}+{{[(n+1)/2]}\choose 1}}{{{[n/2]}\choose 2}+{{[(n+1)/2]}\choose 2}}=
4\frac{n-1+\varepsilon}{n-2+\varepsilon},
\eqno(5)
$$
where $\varepsilon=0$ if $n$ is even and $\varepsilon=1$ if $n$ is odd.\\
Combining (4) with (5), we obtain
$$
5-\frac{m_4}{\alpha_2}\le \alpha_2^{(1)}<4\frac{n-1+\varepsilon}{n-2+\varepsilon}.
\eqno (6)
$$

Denote by $l$ the number of vertices of $P$. 
Denote by $N_4$ the total number 
of vertices of quadrilateral $2$-faces. Clearly, $N_4=4m_4$.
By Lemmas~\ref{quadrilaterals} and~\ref{quad9}
each of $l$ vertices is incident to at most $n+6$ quadrilaterals.
Thus, $N_4\le l(n+6)$ and
we have  $4m_4\le l(n+6)$.
In view of (1) and (5), we have  
$$
\frac{m_4}{\alpha_2}\le \frac{1}{4}\frac{l(n+6)}{\alpha_2}=
 \frac{n+6}{4}\frac{2}{n(n-1)}\alpha_2^{(1)}<\phantom{asasdadsadsadssadfgfgh}
$$
$$
\phantom{asasdadsafghfh}
<\frac{n+6}{2n(n-1)}\frac{4(n-1+\varepsilon)}{(n-2+\varepsilon)}=
2\frac{n+6}{n(n-1)}\frac{(n-1+\varepsilon)}{(n-2+\varepsilon)}.
\eqno (7)
$$

Combining (6) and (7), we obtain
$$
5-\frac{4(n-1+\varepsilon)}{(n-2+\varepsilon)}<\frac{m_4}{\alpha_2}<
2\frac{n+6}{n(n-1)}\frac{(n-1+\varepsilon)}{(n-2+\varepsilon)}.
$$
This implies
$$
(n-6+\varepsilon)n(n-1)<2(n+6)(n-1+\varepsilon).
$$
This is equivalent to $n^2-8n-12<0$ if $n$ is even and
to $n^2-8n-7<0$ if $n$ is odd.
The first inequality has no solutions for $n\ge 10$, and the second one has no solutions for $n\ge 9$.
So, the theorem is proved.\\
{}\hfill\qed

\vspace{25pt}
\noindent
Independent University of Moscow, Russia \\
e-mail: \phantom{ow} felikson@mccme.ru\qquad pasha@mccme.ru

\end{document}

%% file: cox_e.txt
\begin{table}[!h]
\caption{ Connected
elliptic and parabolic Coxeter diagrams are listed in the left and right columns respectively.}
\label{el-par}
\medskip
\begin{center}
\begin{tabular}{|cc@{\quad}|cc|}
\hline
\raisebox{0pt}{${ A_n}$ $(n\ge 1)$}  & 
\raisebox{0pt}{\epsfig{file=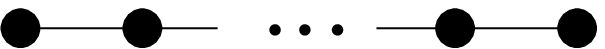,width=0.2\linewidth}}&
\multicolumn{2}{c|}{
\begin{tabular}{cc}
\raisebox{0pt}[20pt][5pt]{${ \widetilde A_1}$} & 
\raisebox{0pt}[20pt][5pt]{\epsfig{file=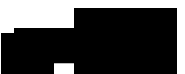,width=0.052\linewidth}}\\
\raisebox{3pt}[15pt][14pt]{${ \widetilde A_n}$ $(n\ge 2)$}  & 
\raisebox{-8pt}[25pt][7pt]{\epsfig{file=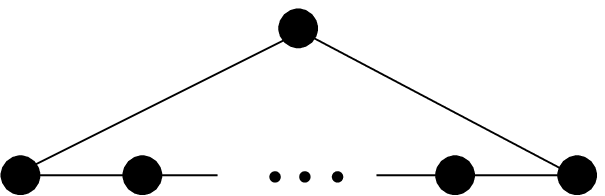,width=0.2\linewidth}}
\end{tabular}
}\\
\hline
\raisebox{-1pt}[23pt][7pt]{${ B_n=C_n}$} & 
\raisebox{-7pt}[23pt][7pt]{\epsfig{file=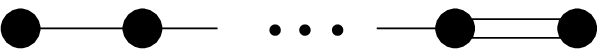,width=0.2\linewidth}}&
\raisebox{7pt}[23pt][7pt]{${ \widetilde B_n}$ $(n\ge 3)$}  & 
\raisebox{-0pt}[30pt][7pt]{\epsfig{file=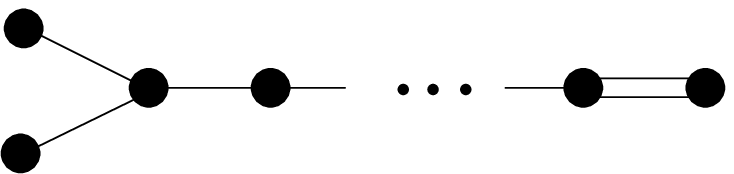,width=0.2\linewidth}}\\
\cline{3-4}
\raisebox{7pt}[15pt][7pt]{ $(n\ge 2)$} & 
&
\raisebox{-0pt}[15pt][7pt]{${ \widetilde C_n}$ $(n\ge 2)$}  & 
\raisebox{-0pt}[15pt][7pt]{\epsfig{file=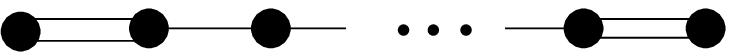,width=0.2\linewidth}}\\
\hline
\raisebox{7pt}[23pt][7pt]{${ D_n}$ $(n\ge 4)$} & 
\raisebox{0pt}[30pt][7pt]{\epsfig{file=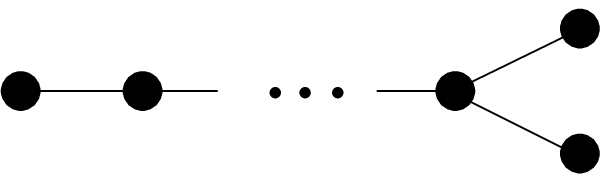,width=0.2\linewidth}}&
\raisebox{7pt}[23pt][7pt]{${ \widetilde D_n}$ $(n\ge 4)$} & 
\raisebox{0pt}[30pt][7pt]{\epsfig{file=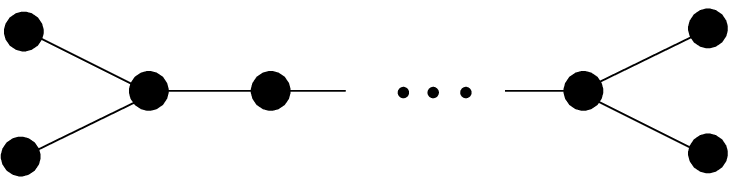,width=0.2\linewidth}}\\
\hline
\raisebox{0pt}[15pt][7pt]{${ G_2^{(m)}}$}  & 
\psfrag{m}{{\scriptsize $m$}}
\raisebox{0pt}[15pt][7pt]{\epsfig{file=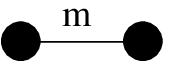,width=0.052\linewidth}}& 
\raisebox{0pt}[15pt][7pt]{${ \widetilde G_2}$} & 
\raisebox{0pt}[15pt][7pt]{\epsfig{file=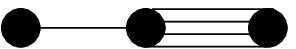,width=0.082\linewidth}}\\
\hline
\raisebox{0pt}[15pt][7pt]{${ F_4}$}  & 
\raisebox{0pt}[15pt][7pt]{\epsfig{file=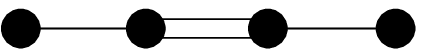,width=0.13\linewidth}}&
\raisebox{-0pt}[15pt][7pt]{${ \widetilde F_4}$} & 
\raisebox{-0pt}[15pt][7pt]{\epsfig{file=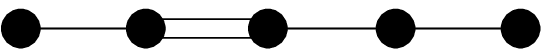,width=0.16\linewidth}}\\
\hline
\raisebox{8pt}[15pt][7pt]{${ E_6}$}  & 
\raisebox{0pt}[30pt][7pt]{\epsfig{file=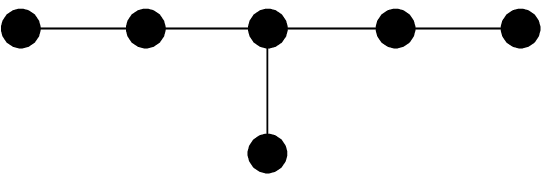,width=0.16\linewidth}}&
\raisebox{8pt}[35pt][7pt]{${ \widetilde E_6}$} & 
\raisebox{-8pt}[40pt][17pt]{\epsfig{file=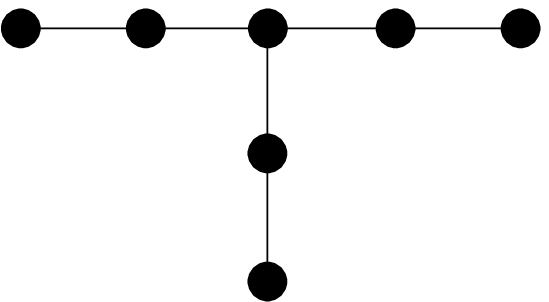,width=0.16\linewidth}}\\
\hline
\raisebox{5pt}[25pt][7pt]{${ E_7}$}  & 
\raisebox{-0pt}[30pt][7pt]{\epsfig{file=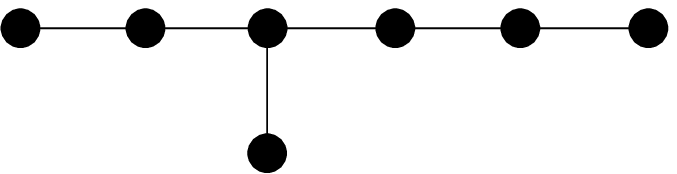,width=0.2\linewidth}}&
\raisebox{5pt}[25pt][7pt]{${ \widetilde E_7}$} & 
\raisebox{-0pt}[25pt][7pt]{\epsfig{file=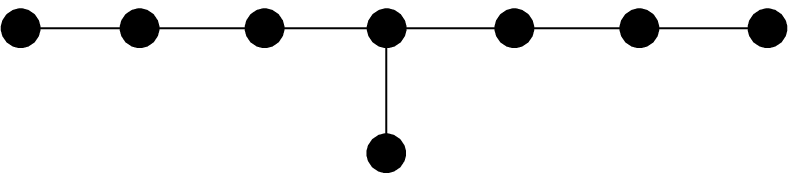,width=0.22\linewidth}}\\
\hline
\raisebox{5pt}[25pt][7pt]{${ E_8}$}  & 
\raisebox{-0pt}[30pt][7pt]{\epsfig{file=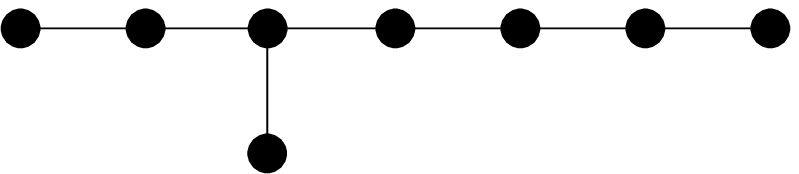,width=0.24\linewidth}}&
\raisebox{5pt}[25pt][7pt]{${ \widetilde E_8}$} & 
\raisebox{-0pt}[25pt][7pt]{\epsfig{file=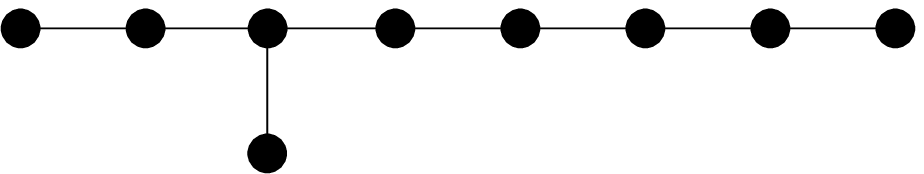,width=0.25\linewidth}}\\
\hline
\raisebox{0pt}[15pt][7pt]{${ H_3}$}  & 
\raisebox{0pt}[15pt][7pt]{\epsfig{file=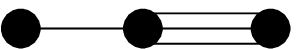,width=0.1\linewidth}}&
& 
\\
\cline{1-2}
\raisebox{0pt}[15pt][7pt]{${ H_4}$}  & 
\raisebox{0pt}[15pt][7pt]{\epsfig{file=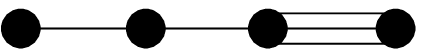,width=0.13\linewidth}}&
& 
\\
\hline
\end{tabular}
\end{center}
\end{table}